\documentclass[a4paper,12pt%
]{amsart}

\usepackage{amsmath, amssymb, amsthm, amscd}
\usepackage{graphicx}
\usepackage{enumerate}
\usepackage{caption}
\usepackage{subcaption}
\usepackage[all,cmtip]{xy}

\usepackage{algorithm}% http://ctan.org/pkg/algorithms
\usepackage{algpseudocode}% http://ctan.org/pkg/algorithmicx

\usepackage[english]{babel}
\usepackage[latin1]{inputenc}
\usepackage[all,cmtip]{xy}
\usepackage{times}
\usepackage{cancel}

\newcounter{theorem}

\newtheorem{defn}[theorem]{Definition}

\theoremstyle{remark}
\newtheorem*{remark*}{Remark}

\newtheorem{rle}[theorem]{Rule}

\numberwithin{equation}{section}
\numberwithin{theorem}{section}

\newcommand{\defemph}{\emph}

\newcommand{\R}{\mathbb{R}}

\newcommand{\PT}{R} %  a set of prototiles
\newcommand{\Rhombs}[1]{\PT(#1)} %  prototiles with special angles
\newcommand{\Rhomb}[1]{r(#1)} %  rhomb with angle i pi /n
\newcommand{\tot}[1]{\sum #1} %  the sum of the direction vectors of an edge sequence
\newcommand{\Patch}{P} %  a patch
\newcommand{\tile}{t} %  a single tile
\newcommand{\prototile}{p} %  a single prototile
\newcommand{\sub}{\varphi} %  a substitution
\newcommand{\infl}{\lambda} %  an inflation factor
\newcommand{\edge}{e} %  the edge of a tile
\newcommand{\seq}{\sigma} %  an edge sequence 
\newcommand{\pt}{x} % a point in the plane
\newcommand{\PatchSet}{\mathcal{P}} %  a set of patches

\title[Rhomb tilings]{On substitution tilings of the plane with $n$-fold rotational symmetry}

\author{Gregory R. Maloney}
\address{Newcastle University}
\thanks{}

\begin{document}

\subjclass[2010]{Primary: 52C20, 05B45, 52C23.  
Secondary: 52C30}
\keywords{Rhombs, Tiling, Algorithm}
\date{\today}

\begin{abstract}
A method is described for constructing, with computer assistance, planar substitution tilings that have $n$-fold rotational symmetry.  
This method uses as prototiles the set of rhombs with angles that are integer multiples of $\pi/n$, and includes various special cases that have already been constructed by hand for low values of $n$.  
An example constructed by this method for $n = 11$ is exhibited; this is the first substitution tiling with $11$-fold symmetry appearing in the literature.  
\end{abstract}

\maketitle

\section{Introduction}\label{SEC:intro}

Rotational symmetry is one of the most distinctive qualities that aperiodic planar tilings can have.  
The most famous aperiodic family of tilings is the family of Penrose tilings \cite{P:penrose}, which possess fivefold rotational symmetry.  
The Ammann-Beenker tilings \cite{B:ammann-beenker} are arguably the second-most famous aperiodic tilings, and they possess eightfold rotational symmetry.  
Aperiodic tilings are often used as models of quasicrystals; when Dan Shechtman and his co-authors made the Nobel prize-winning discovery of quasicrystals in \cite{S:nobel}, the first sign by which they knew that they had found something special was the presence of tenfold rotational symmetry in the X-ray diffraction patterns.  

Yet there is a lack of known examples of aperiodic planar tiling families with higher orders of rotational symmetry (see the Tilings Encyclopedia \cite{TE} for an extensive list of substitution tilings, including most of the ones cited here).  
It is true that there are projection tilings possessing arbitrary orders of rotational symmetry---see \cite{BKSZ:projection} for a method with $n = 5$ that generalizes easily to arbitrary $n$---but the same cannot be said for substitution tilings.  
And under a certain assumption substitution tilings possess the desirable property of \defemph{repetitivity}, which implies in particular that any local configuration of tiles with rotational symmetry that appears in one such tiling must appear in every ball of sufficiently large radius in every such tiling.  
Therefore local patterns with rotational symmetry must appear quite regularly in substitution tilings if they appear at all, whereas in projection-method tilings heuristic arguments and empirical observation suggest that patterns with high orders of rotational symmetry must necessarily be sparse (see \cite[Remark 7.12, p. 300f]{BG:book} and \cite{MSRRSB:lasers}).  

The purpose of this work is to provide a method for finding substitution tilings with high orders of rotational symmetry.  
In practical terms, ``high orders'' of rotational symmetry are numbers in the range eleven to nineteen; indeed, to my knowledge, the highest known order of rotational symmetry of any substitution tiling is twelve \cite{G:shield, S:socolar, WSI:dodecagonal}, and there are no known substitution tilings in the literature with elevenfold symmetry (but see \cite{R:subrosa} for some unpublished examples similar to the ones presented here, including one with elevenfold symmetry).  
Therefore the substitution tiling space with elevenfold symmetry described in Section \ref{SEC:example} of this paper is the first of its kind to be published.  

Much previous work has been done on this topic, some of which claims to produce planar substitution tiling spaces with $n$-fold rotational symmetry for arbitrary $n$.  
My claim that there are no previously-known examples of substitution tilings with elevenfold rotational symmetry would seem to contradict this previous work, but in fact there is no contradiction; instead the confusion lies in differing interpretations of $n$-fold symmetry in a tiling space.  
Let us therefore be more explicit about what $n$-fold symmetry means in the current context.  

Two different repetitive tilings are locally indistinguishable from one another: any bounded patch that appears in one of them also appears in the other.  
Therefore it is natural to consider the space of all tilings arising from a given substitution rule, rather than just a single such tiling.  
Many results in the literature use the term $n$-fold symmetry to refer to tiling \emph{spaces} that are $n$-fold symmetric; that is, a rotation of any tiling in the space by $\pi/n$ yields another tiling in the space.  
This is not the same as saying that the tiling space contains an individual tiling that is invariant under $n$-fold rotation, which is precisely the property that is sought here.  

A tiling that is invariant under $n$-fold rotation contains bounded rotation-invariant patches of arbitrary radius centred on its centre of rotation.  
Conversely, the existence of a non-trivial rotation invariant patch implies the existence of a rotation-invariant tiling in the following way.  
Place the patch over the origin and apply the substitution rule to it repeatedly until two patches are found, one of which is contained in the other.  
The smaller of these is the seed for a tiling that is a fixed point of some power of the substitution; this tiling will necessarily have rotational symmetry.  

But we must be careful what we mean by a ``rotation-invariant patch.''  
It is a common situation to have two congruent tiles that are distinguished from each other by giving them different labels.  
Thus a rotation-invariant patch is collection of tiles for which there is a rotation that sends each tile to another tile in the collection that not only is congruent to the original tile, but also has the same label.  
Since the substitution behaves in the same way for two tiles with the same label, this guarantees that the image of the patch after several substitutions will also be rotation-invariant.  
This condition is not met in \cite[Remark 6.3]{DN:n-fold}, in which there are patches consisting of fourteen isosceles triangles that appear to have fourteenfold rotational symmetry, but that in fact only have twofold symmetry.  
This is because some of the tiles are reflected copies of the others, and hence they have different labels.  

The goal of this work is to find substitution tilings that contain patches of bounded size that are invariant under $n$-fold rotation in the sense just described.  
There are two types of tiles that are naturally suited to this purpose.  
The first is the collection of triangles with angles that are integer multiples of $\pi/n$; the second is the set of rhombs with the same angles.  
These tiles have been used in substitution tilings before (see \cite{DN:n-fold, F:thesis, G:multi, G:fivefold, GKM:search} for triangles and \cite{P:penrose, B:ammann-beenker, H:rhomb, R:subrosa} for rhombs), but except for a few specific examples with low values of $n$, it is always the tiling space itself that exhibits rotational symmetry, and not any individual tiling.  
The examples produced here do have $n$-fold symmetry; specifically, they contain patches consisting of stars of $n$ congruent rhombs, each of which has angles $2\pi/n$ and $(n-2)\pi/n$.  

Most previous work on this subject has focused on examples that can be created by hand.  
It is possible to create substitutions with fivefold or sevenfold symmetry by hand, putting families of tiles together to form larger copies of themselves, but for higher orders of symmetry the scale becomes too large for this method to work.  
A common approach to this problem of scale has been to take a known substitution rule with fivefold or sevenfold symmetry and then to generalize it to higher $n$ in some way.  
This has led to the discovery of various infinite families of tiling spaces, but in all cases the desirable property of having an individual tiling with $n$-fold symmetry is lost upon passage to higher $n$.  

The approach here is not to use any specific fivefold or sevenfold example as a model, but rather to apply in a systematic way the ad hoc methods that have been used to discover such examples.  
The scale of this problem at higher values of $n$ means that it is necessary to use a computer.  
In particular there is an algorithm of Kannan, Soroker, and Kenyon \cite{KS:parallelograms,K:parallelograms} that will produce a substitution rule of a given size and shape, if it is possible to do so.  
Then a result of Kenyon says that any other substitution rule with that size and shape can be obtained from this one by a series of transformations, called rotations.  
Again with the help of a computer we can transform the initial substitution rule via rotations until it produces $n$-fold symmetric patches.  
This is not possible for every size and shape of substitution, but for large enough substitutions it works.  

It should be mentioned that the method described here relies mostly on existing algorithms and ideas.  
The chief innovation is Rule \ref{RULE} in Section \ref{SEC:definitions}.  
Nevertheless, it should not be dismissed just because its ingredients are not novel.  
Indeed, these ingredients have been known for many years, during which many researchers have tried without success to produce $n$-fold symmetric substitution tilings for large $n$; the fact that this has been achieved here is evidence that the combination of these ingredients in this particular way is a significant development.  

\section{Definitions and notation}\label{SEC:definitions}

Let us suppose henceforth that an odd integer $n\geq 3$ is fixed.  
The case when $n$ is even is addressed in Section \ref{SEC:even}.  

\begin{defn}\label{DEF:tiling}
A \defemph{tile} is a subset of $\R^d$ homeomorphic to the closed unit disk.  
A \defemph{patch} is a collection of tiles, any two of which intersect only in their boundaries.  
The \defemph{support} of a patch is the union of the tiles that it contains.  
A \defemph{tiling} is a patch, the support of which is all of $\R^d$.  
\end{defn}

\begin{figure}
\includegraphics[width=\textwidth]{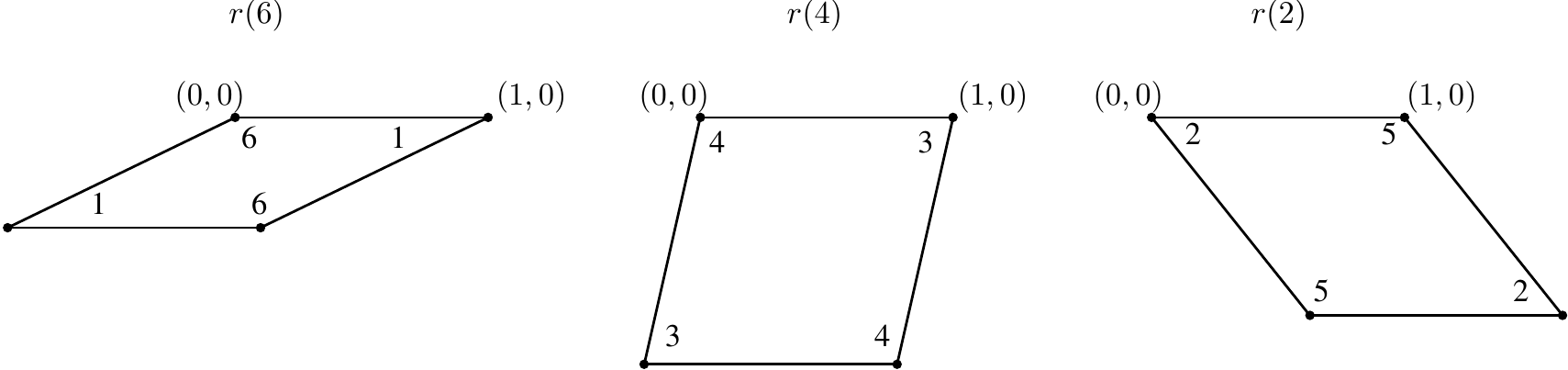}
\caption{\label{FIG:prototiles}The prototiles $\Rhombs{7}$.  Integers in the corners represent angle measures expressed as integer multiples of $\pi/7$.}
\end{figure}

Let us restrict our attention to the case $d = 2$, and let us consider a restricted set of tiles consisting of all rhombs with sides of unit length and angles that are integer multiples of $\pi/n$.  
Let us select a representative for each isometry-equivalence class of such tiles: for each even integer $i < n$, let $\Rhomb{i}$ denote the rhomb with vertices $(0,0)$, $(1,0)$, $(1-\cos (i\pi/n),-\sin (i\pi/n))$, and $(-\cos (i\pi/n),-\sin (i\pi/n))$.  
Let $\Rhombs{n} = \{ \Rhomb{i} \ | \ 1\leq i \leq n/2\}$.  
Then $\Rhombs{n}$ contains $\lfloor \frac{n}{2}\rfloor$ elements if $n$ is odd, where $\lfloor \cdot \rfloor$ denotes the floor function.  
The elements of $\Rhombs{n}$ are called \defemph{prototiles}, and we denote by $\PatchSet(\Rhombs{n})$ the set of all patches consisting entirely of tiles congruent to these prototiles.  
Figure \ref{FIG:prototiles} depicts the rhombs $\Rhombs{7}$.  

This work deals with tilings that arise from substitutions.  
A substitution is traditionally defined to be a map on a set of prototiles that assigns to each prototile $\prototile$ a patch, the union of which is $\infl\prototile$ for some $\infl>1$, called the \defemph{inflation factor}.  
However, in the present context we will need to relax this definition slightly; in particular, we will need to deform the edges of $\infl\prototile$.  
For this it will be necessary to introduce some new notation and terminology.  

For an integer $k$, let $\edge_k$ denote the direction vector beginning at $(0,0)$ and ending at $(\cos (k\pi/n),\sin (k\pi/n))$.

\begin{defn}\label{DEF:edge-sequence}
An \defemph{edge sequence} is a finite sequence of integers $\seq = (k_i)_{i=1}^m$.  
A \defemph{standard} edge sequence is one for which $\sum_{i=1}^m k_i = 0$ and $|k_i| < n/2$ for all $1\leq i\leq m$.  
To an edge sequence $\seq$ and a point $\pt\in\R^2$ we associate a sequence $E_\seq(\pt)$ of line segments, each of which begins at the end of the previous one and has direction vector $\edge_{k_j}$.  
The first line segment of $E_\seq$ begins at the point $\pt$.  

Given an edge sequence $\seq$, let $\tot{\seq}$ denote the point $\sum_{i=1}^m e_{k_i}$.  
Then the endpoint of $E_\seq(0)$ is $\tot{\seq}$.  

Given an integer $j$, define the \defemph{rotation of $\seq$ by $j$} to be the edge sequence $\seq(j) := (k_i+j)_{i=1}^m$.  
\end{defn}

The conditions in the definition of a standard edge sequence $\seq$ imply that $\tot{\seq}$ lies on the positive $x$-axis.  

Now we can use a standard edge sequence $\seq$ to distort the boundaries of the rhombs in $\Rhombs{n}$.  

\begin{defn}\label{DEF:boundary}
Let $\seq$ be a standard edge sequence $\seq$ and $i < n$ an even integer.  
The \defemph{$\seq$-boundary of $\Rhomb{i}$} is $\seq \Rhomb{i} := E_\seq(0) \cup E_\seq(\tot{\seq(-i)})\cup E_{\seq(-i)}(0)\cup E_{\seq(-i)}(\tot{\seq})$.  
\end{defn}

\begin{figure}
\includegraphics[width=\textwidth]{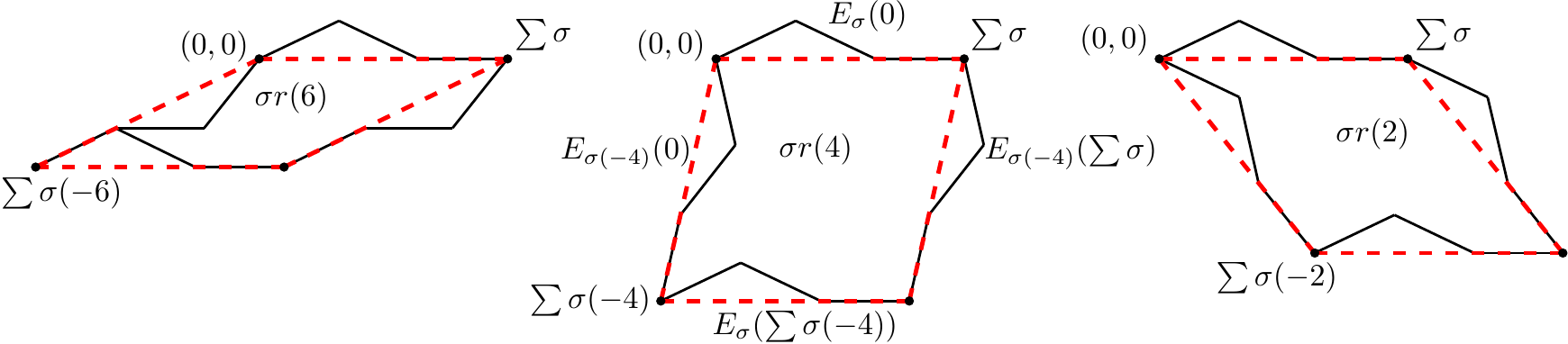}
\caption{\label{FIG:boundaries}$\seq$-boundaries of $\Rhombs{7}$ for $\seq = (1,-1,0)$.   The outlines of the inflated copies of these rhombs are depicted as red dashed lines.  The edges of $\seq \Rhomb{4}$ are labelled.}
\end{figure}

The end points of the various line segment sequences in the $\seq$-boundary of $\Rhomb{i}$ form the vertices of an inflated copy of $\Rhomb{i}$.  
Figure \ref{FIG:boundaries} depicts the $(1,-1,0)$-boundaries of $\Rhomb{6}$, $\Rhomb{4}$, and $\Rhomb{2}$ for $n = 7$.  
The outlines of the inflated copies of these rhombs are depicted as dashed lines.  

The $\seq$-boundary of $\Rhomb{i}$ is the image of a closed curve, specifically the curve defined by traversing at a constant speed $E_\seq(0)$ and $E_{\seq(-i)}(\tot{\seq})$ in the positive direction, then $E_\seq(\tot{\seq(-i)})$ and $E_{\seq(-i)}(0)$ in the negative direction.  
This curve is not necessarily simple---see for example $\seq \Rhomb{6}$ in Figure \ref{FIG:boundaries}, which backtracks on itself near the bottom-left corner.  

The goal here is to define a substitution rule on $\Rhombs{n}$ that sends each prototile $\Rhomb{i}$ to a patch, the support of which is the closure of the planar region inside of $\seq \Rhomb{i}$.  
This is only possible if it is understood what is meant by ``the planar region inside of $\Rhomb{i}$.''  
When $\seq\Rhomb{i}$ is a simple curve, the meaning of this is clear.  
It is also clear that $\seq \Rhomb{6}$ in Figure \ref{FIG:boundaries}, which involves some backtracking, has a well-defined region inside of it.  
The notion of an inside makes sense for a more general class of curves, although they are somewhat difficult to describe.  

\begin{defn}\label{DEF:good-boundary}
Let $E$ be a closed path of directed line segments of unit length.  
$E$ is a \defemph{good curve} if it satisfies the following.  
\begin{enumerate}
\item  If any two line segments overlap at a point other than their endpoints, then both segments are the same with opposite orientations; and
\item  it is possible to perform small homotopies to separate all pairs of overlapping oppositely-oriented line segments to obtain a planar-embedded graph for which, at any vertex where four edges meet, those four edges are alternating in- and out-edges of that vertex.  
\end{enumerate}
\end{defn}

\begin{figure}
\includegraphics[width=\textwidth]{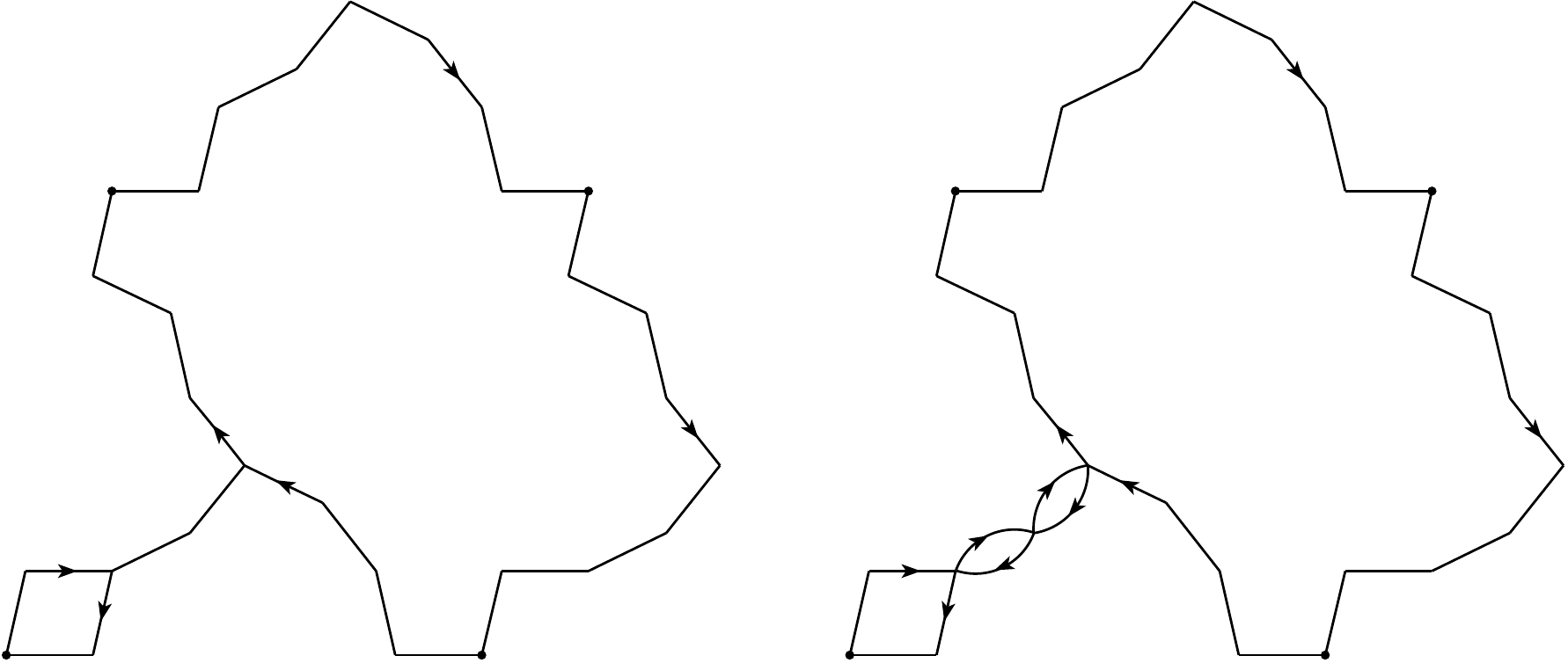}
\caption{\label{FIG:strange-boundary}The $(0,3,1,2,-1,-2,-3,0)$-boundary of $\Rhomb{4}$.  The picture on the right shows the images of the double lines under homotopies.  }
\end{figure}

The curve in Figure \ref{FIG:strange-boundary} is a good curve, as can be seen on the right side of the figure, in which homotopies fulfilling the requirements of Definition \ref{DEF:good-boundary} are exhibited.  

\begin{defn}\label{DEF:substitution}
A \defemph{substitution} on $\Rhombs{n}$ is a map $\sub : \Rhombs{n} \to \PatchSet(\Rhombs{n})$ for which there is an edge sequence $\seq$ such that each $\seq \Rhomb{i}$ is a good curve, the closure of the inside of which is the support of $\sub (\Rhomb{i})$.  
\end{defn}

\begin{figure}
\includegraphics[width=\textwidth]{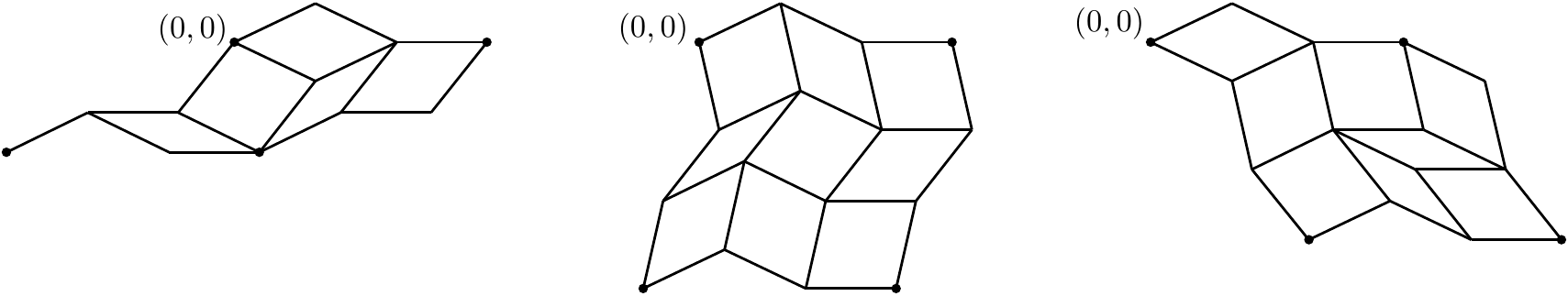}
\caption{\label{FIG:substitution}A substitution using the $(1,-1,0)$-boundaries of $\Rhombs{7}$.}
\end{figure}

Figure \ref{FIG:substitution} depicts a substitution on $\Rhombs{7}$ that uses the $(1,-1,0)$-boundaries; these boundaries appeared earlier in Figure \ref{FIG:boundaries}.  

The patches in Figure \ref{FIG:substitution} are ambiguous because they do not specify how the tiles are obtained as images of the prototiles under isometries.  
The prototiles in $\Rhombs{7}$ are all self-symmetric under reflection and rotation by a half turn.  
Therefore each rhomb that appears in the figure could be the image of a prototile under any one of four isometries, two of which are orientation preserving and two of which are orientation reversing.  
In order to determine a substitution rule fully, we must specify which of these isometries to use for each tile in each patch.  

The rule for specifying these isometries is simple but important; without this rule it would not be possible to iterate $\sub$.  
Let us assign an orientation to each line segment $L$ in the plane that is parallel to one of the direction vectors $\edge_k$.  

\begin{rle}\label{RULE}
If the $L$ is parallel to $\edge_k$ for an even integer $0\leq k < n$, then give $L$ the same orientation as $k$; otherwise give it the opposite orientation.  
\end{rle}

This gives orientations to the edges of all of the rhombs in $\Rhombs{n}$, and also to the rhombs in the patches $\sub (\Rhomb{i})$.  
These orientations have the property that, if two line segments have a common end point, then if the angle between them is an even multiple of $\pi/n$ the line segments are oriented either both toward or both away from their common end point; if the angle between them is an odd multiple of $\pi/n$ then one of the line segments points toward the common end point and the other one points away from it.  

\begin{figure}
\includegraphics[width=\textwidth]{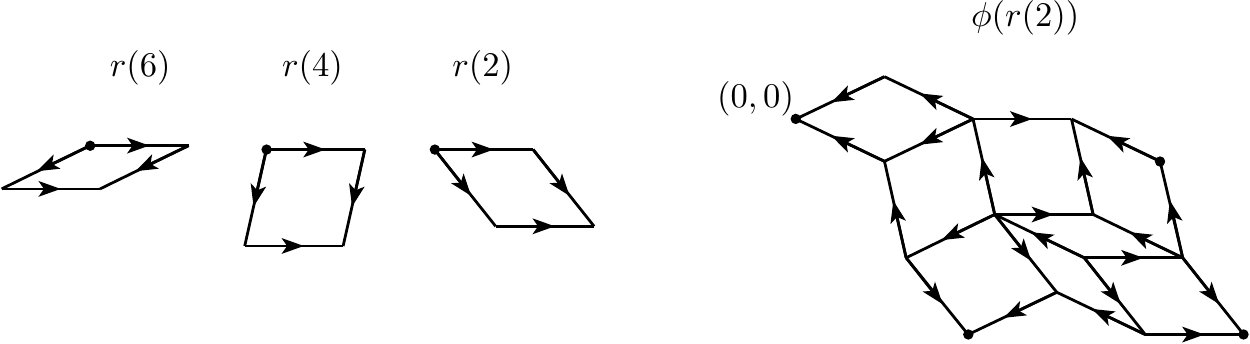}
\caption{\label{FIG:orientations}The prototiles $\Rhombs{7}$ and the substituted image $\sub(\Rhomb{2})$, all with their edge orientations.  }
\end{figure}

These orientations determine, for each tile in a patch, a unique orientation-preserving isometry that carries that tile to its associated prototile.  
Figure \ref{FIG:orientations} depicts the three prototiles in $\Rhombs{7}$ and the patch $\sub(\Rhomb{2})$, all with their edge orientations.  

The reason for establishing these orientation rules is that it allows us to extend $\sub$ to arbitrary patches, which in turn means that we can iterate $\sub$.  
This is a crucial procedure in the construction of substitution tilings (see \cite[Chapter 6]{BG:book}, where the less standard but more precise term \defemph{inflation tilings} is used).  

First let us extend the domain of a substitution $\sub$ to include all tiles of the form $g(\Rhomb{i})+v$, where $g$ is an orthogonal transformation of the plane and $v$ is a translation vector.  
For this it will be necessary to generalize the notion of an \defemph{inflation factor}.  

\begin{defn}\label{DEF:inflation-factor}
Let $\sub$ be a substitution with associated boundary sequence $\seq$.  
Then the \defemph{inflation factor} of $\sub$ is the positive number $\infl = \| \tot{\seq}\| $.  
\end{defn}

Then given a substitution $\sub$ with inflation factor $\infl$ and a transformed prototile $g(\Rhomb{i})+v$, let us define 
\begin{align*}
\sub(g(\Rhomb{i})+v) := \{ g(\tile) + \infl v \ | \ \tile \in \sub(\Rhomb{i})\}.
\end{align*}
Finally, for a patch $\Patch$ consisting of rhombs with edges parallel to the direction vectors $\edge_k$, let us define 
\begin{align*}
\sub(\Patch) := \{ \sub(\tile) \ | \ \tile \in \Patch\}.
\end{align*}

The $\seq$-boundaries of the prototiles in $\Rhombs{n}$ are consistent with the orientations stipulated by Rule \ref{RULE}.  
This means that the image under $\sub$ of the patch $\Patch$ is also a patch---boundaries of two tiles that share an edge are distorted in the same way, so there is no overlap (of positive measure).  

In particular every set $\sub(\Rhomb{i})$ is a patch consisting of rhombs with edges parallel to the direction vectors $\edge_k$, so we can define $\sub^2(\Rhomb{i}) := \sub(\sub(\Rhomb{i}))$ and so on.  

% Tiling of the plane.
% Substitution, substitution tiling.
% Tiling space?
% Our particular prototile set.  
% Notation for substitution.  

\section{The method}\label{SEC:method}

\subsection{The Kannan-Soroker-Kenyon criterion}\label{SUBSEC:criterion}

To find a substituion on the prototiles $\Rhombs{n}$ that has $n$-fold rotational symmetry, the first step is to choose a standard edge sequence $\seq$.  
This involves some work; let us discuss this further in Section \ref{SEC:boundaries}.  

The next step is to determine if there is at least one substitution $\sub$ associated to $\seq$.  
This part also requires work, but fortunately Kannan and Soroker \cite{KS:parallelograms} and Kenyon \cite{K:parallelograms} have developed a simple criterion for determining if such a substitution exists; this criterion can be checked with a computer.  
More specifically, the Kannan-Soroker-Kenyon (KSK) criterion is a test that determines if the region bounded by a curve consisting of unit-length line segments is tilable by the prototiles of $\Rhombs{n}$; if each boundary $\seq\Rhomb{i}$ satisfies this criterion, then a substitution exists.  

The KSK criterion can be applied to any polygonal boundary, but in the special case of the boundaries $\seq \Rhomb{i}$ it can be expressed in particularly simple terms.  
This involves \defemph{pseudoline arrangements} (see \cite{FH:worms}, where they are called \defemph{worms} following Conway).  

If the region inside of $\seq\Rhomb{i}$ can be tiled by rhombs from $\Rhombs{n}$, then it has an associated pseudoline arrangement.  
A pseudoline is a smooth planar curve, and a pseudoline arragement is a collection of pseudolines, any two of which, if they intersect, do so at most once, and cross each other at their point of intersection.  
Let us impose a further condition in the definition of a pseudoline arrangement that is not standard in the literature: no three pseudolines may share a common intersection point.  

Then a pseudoline arrangement is associated to a tiling of the inside of the region bounded by $\seq\Rhomb{i}$ in the following way.  
Group the rhombs in the tiling into inclusion-maximal lines of linked rhombs with a common edge direction, each of which shares a parallel edge with the next.  
Every rhomb falls into exactly two such families.  
A pseudoline arrangement is obtained by drawing a pseudoline through each such family of rhombs.  

\begin{figure}
\includegraphics[height=4.5cm]{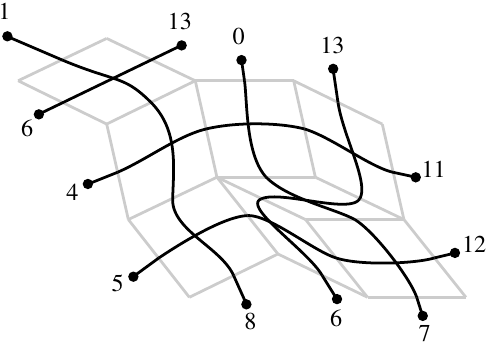}
\includegraphics[height=4.5cm]{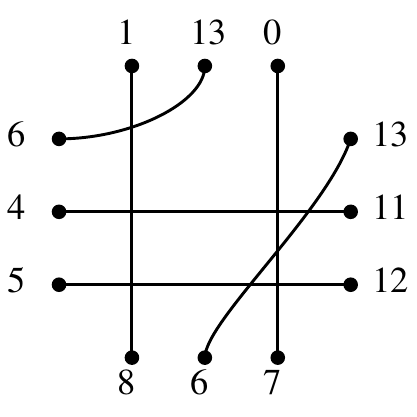}
\caption{\label{FIG:pseudoline}Pseudoline arrangement for $\sub(\Rhomb{2})$ from Figure \ref{FIG:substitution}.  The picture on the right is the same arrangement with the endpoints of the pseudolines shifted to lie on a square.}
\end{figure}

Figure \ref{FIG:pseudoline} depicts such a pseudoline arrangement for the patch $\sub(\Rhomb{2})$ from Figure \ref{FIG:substitution}.  
These pseudoline arrangements are considered equivalent up to diffeomorphisms that preserve the ``no triple intersections'' property, so the pseudoline arrangement on the left of that figure is equivalent to the one on the right, which has been deformed for clarity of drawing, and from which the background tiling has been removed.  

The drawing on the right of Figure \ref{FIG:pseudoline} also has numbers at the ends of the pseudolines indicating the angles that the corresponding line segments of $\seq\Rhomb{2}$ make with the positive $x$-axis.  
Each pseudoline joins two edges angles that differ by $\pi$, so labels of the end points of a single pseudoline always sum to $2n$.  
This is true for any $n$ and for any pseudoline arrangement arising from a tiling of the inside of $\seq\Rhomb{i}$.  

Moreover, Definition \ref{DEF:edge-sequence} of a standard edge sequence implies that each angle $0\leq j < 2n$ will only occur as a label of line segments on two of the four sides of $\seq\Rhomb{i}$, and those two sides cannot be opposite.  
Therefore, when reading the labels cyclically around the edge of $\seq\Rhomb{i}$, the labels $j$ will not appear interspersed with the labels $n+j$.  
Also, two pseudolines with ends labelled $j$ and $n+j$ cannot intersect one another in an arrangement arising from a tiling.  
These two facts together imply that there is a unique way to partition the line segments in $\seq\Rhomb{i}$ into pairs connected by a common pseudoline.  

Whether a tiling of the inside of $\seq\Rhomb{i}$ exists or not, it is possible to label the edge segments with their angles modulo $2n$, and to determine which must be connected to which by pseudolines.  
It is also possible to determine which pairs of pseudolines must intersect each other; the only thing that cannot be determined just from looking at the boundary $\seq\Rhomb{i}$ is the order in which the pseudolines intersect.  
But it is not necessary to determine this in order to know whether a tiling of the inside of $\seq\Rhomb{i}$ exists or not.  

If a tiling of the region bounded by $\seq\Rhomb{i}$ exists, then an intersection of two pseudolines in the associated pseudoline arrangement corresponds to a rhomb in the tiling.  
The angles in that rhomb can be calculated by subtracting the angles of the endpoint labels of the two pseudolines: pick an endpoint label of one of the lines and subtract it modulo $2n$ from the next label, in counterclockwise order, of the other line, to obtain one of the angles of the rhomb (expressed as an integer multiple of $\pi/n$).  
The KSK criterion simply says that a tiling of the region bounded by $\seq\Rhomb{i}$ exists if and only if the subtractions for all such pairs of crossing pseudolines yield angles less than $\pi$.  

\begin{figure}
\includegraphics[width=0.45\textwidth]{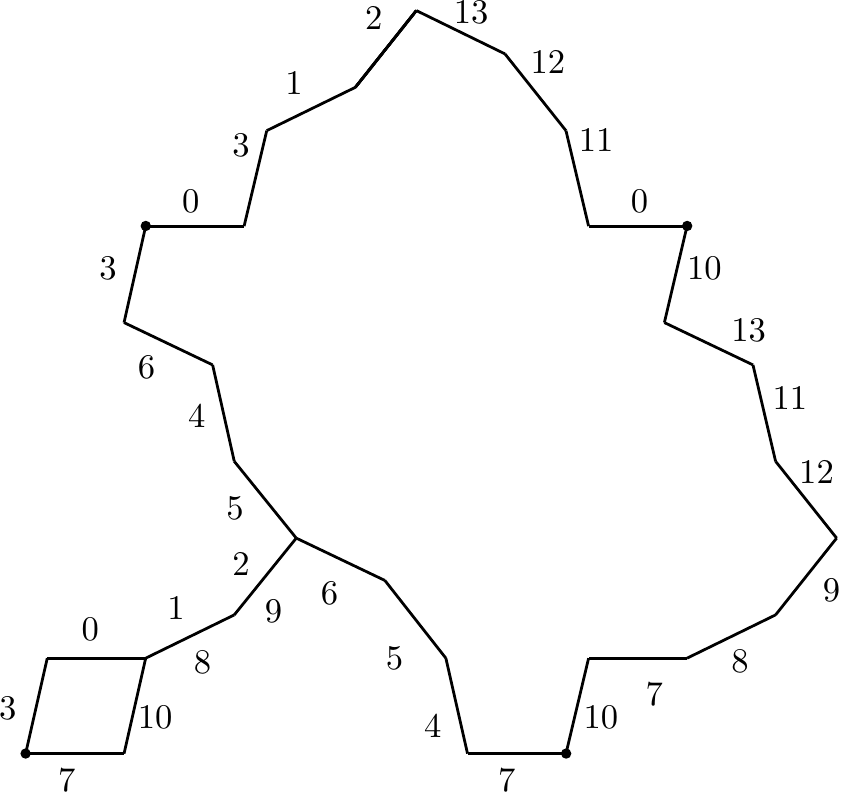}
\includegraphics[width=0.45\textwidth]{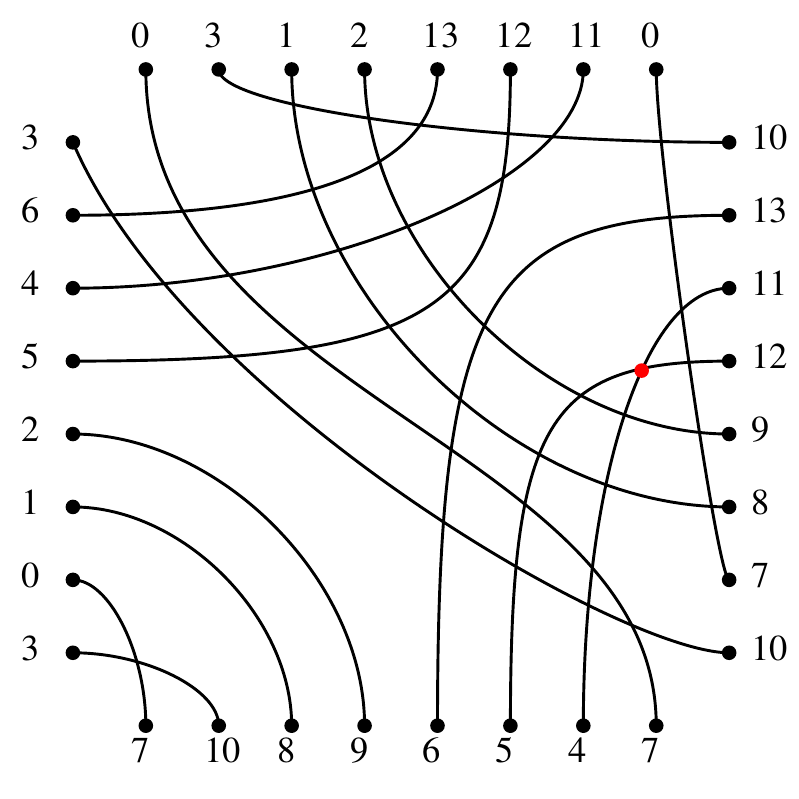}
\caption{\label{FIG:untilable}Pseudoline arrangement for the $(0,3,1,2,$ $-1,-2,-3,0)$-boundary of $\Rhomb{4}$ from Figure \ref{FIG:strange-boundary}.  The red dot indicates the crossing pseudoline pair that violates the KSK criterion, and numbers indicate angles that line segments make with the positive $x$-axis.}
\end{figure}

Figure \ref{FIG:untilable} depicts the $(0,3,1,2,-1,-2,-3,0)$-boundary of $\Rhomb{4}$ from Figure \ref{FIG:strange-boundary} along with an associated pseudoline arrangement.  
The region inside this boundary cannot be tiled by the rhombs in $\Rhombs{7}$, as can be seen from the pseudoline arrangement, which contains a red dot indicating the crossing pseudoline pair that violates the KSK criterion.  
There are many pseudoline arrangements that we could draw for this boundary, and all of them would contain an intersecting pair of pseudolines connecting the same pairs of end points, so none of them is a pseudoline arrangement arising from a tiling by the elements of $\Rhombs{7}$.

\subsection{Rotations}\label{SUBSEC:hex-moves}

Not only did Kannan, Soroker and Kenyon give a criterion to determine if there is a tiling of the inside of $\seq\Rhomb{i}$, they also provided an algorithm to produce such a tiling when it exists.  
By applying this algorithm to each boundary $\seq\Rhomb{i}$, we can produce a substitution.  

\begin{figure}
\includegraphics[width=\textwidth]{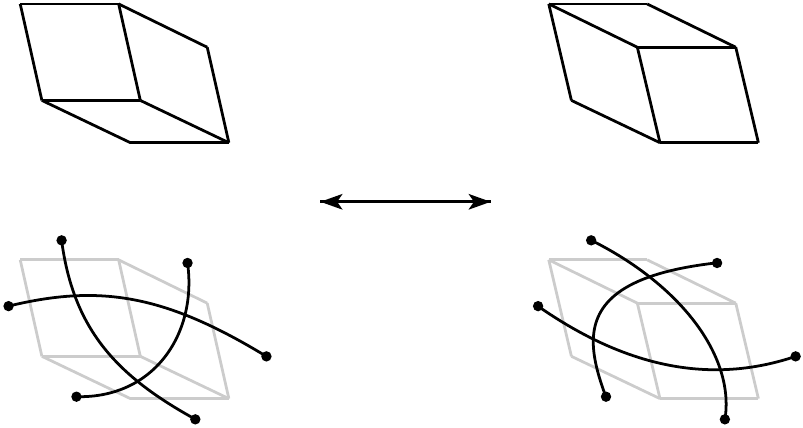}
\caption{\label{FIG:hex-move}A rotation.  The corresponding change in the pseudoline arrangement is depicted underneath.}
\end{figure}

The final step of the method is to modify this substitution to produce a substitution that yields rotational symmetry.  
Kenyon \cite{K:parallelograms} showed that every tiling of the region bounded by $\seq\Rhomb{i}$ using rhombs in $\Rhombs{n}$ is related to every other such tiling by a sequence of simple changes, called rotations.  
Rotations are defined on patches that are triples of rhombs, any two of which share an edge.  
The support of such a triple is a hexagon with opposite edges parallel and of equal length.  
A rotation replaces such a triple with another triple of rhombs each of which is a translate of one of the original three rhombs by a unit vector that is a common edge of the other two rhombs.  
The new triple has the same support as the old one.  
Figure \ref{FIG:hex-move} depicts a rotation, with the corresponding change in the pseudoline arrangement underneath.  

If we represent the tilings $\seq\Rhomb{i}$ using an appropriate data structure, then it is possible to implement rotations with a graphical interface, which makes it easy to find patches with the desired symmetry properties.  
In particular, it is possible to arrange for stars consisting of copies of $\Rhomb{2}$ to appear inside of the patches $\sub(\Rhomb{i})$.  
Another strategy is to put copies of $\Rhomb{2}$ in all the corners of each patch $\sub(\Rhomb{i})$ so that, after a second application of $\sub$, stars appear in the corners of the first supertiles.  
This is what has been done in the example with elevenfold symmetry in Figure \ref{FIG:example-prototiles}.  

Not all edge sequences yield substitutions, and not all edge sequences that yield substitutions yield substitutions that produce $n$-fold rotational symmetry.  
Nevertheless, by passing to large enough edge sequences, it seems that it is always possible to find the desired symmetry.  

\subsection{Variants}\label{SUBSEC:variants}

Let us note here that some well-known substitution tilings arise from very similar methods, but with some important differences.  

The Penrose rhomb substitution \cite{P:penrose} uses the prototile set $\Rhombs{5}$ and also involves tiling inflated prototiles with distorted boundaries, but pairs of opposite edges do not share the same distortion.  
Indeed, for each pair of opposite edges in a prototile, one edge is distorted with the edge sequence $(1,-1)$ and the other is distorted with $(2,-2,0)$.  
This exploits the fact that, if $\seq_e$ and $\seq_o$ are sequences of all the even and all the odd numbers, respectively, between $-n/2$ and $n/2$, then the vectors $\tot{\seq_e}$ and $\tot{\seq_o}$ are the same.  
This follows easily from the fact that the sum of the edge vectors in a regular $n$-gon is the zero vector.  

This approach works for $\Rhombs{5}$, but it is very difficult to generalize because if opposite edges are not distorted in the same way, then Rule \ref{RULE} is no longer sufficient to guarantee that the image of an arbitrary  patch under $\sub$ is still a patch of non-overlapping tiles.  
Instead one must check very carefully to be sure that any two rhombs that meet in $\sub(\Rhomb{i})$ have the same distortion on their common edge.  

Another well-known substitution that arises from a similar method is the binary substitution \cite{LB:binary}, which also uses the prototiles $\Rhombs{5}$ with distorted boundaries.  
The difference in this case is that the boundaries are distorted using sequences of unit-length line segments that are parallel to the edges of a regular decagon, not just a regular pentagon.  
In other words, the edge sequence uses not only integers, but also half integers.  
By stretching Definition \ref{DEF:edge-sequence}, one could say that the binary tiling uses the edge sequence $(1/2,-1/2)$.  

Another class of examples that use half-integer edge sequences appears in \cite{R:subrosa}.  
These examples use a family of edge sequences parametrized by $n$, and, at least for values of $n$ up to and including $n = 11$, they yield boundaries that satisfy the KSK criterion, and therefore give rise to substitution rules.  
Moreover, these substitution rules exhibit $n$-fold rotational symmetry, at least up to $n = 11$.

% Distorted rhomb boundaries.
% Oriented edges, even and odd.  
% If we can find tilings of all distorted inflated rhombs by other rhombs, then we have a substitution.  
% Kannen-Sorokar-Kenyon algorithm (check spellings of names) tells us if a region is tilable based on its boundary.
% Given an edge substitution, apply the algorithm to all n-fold rhombs.  
% If they're all tilable, you've got a substitution; otherwise you don't.  
% Once you know a family of distorted (deformed?) rhombs is tilable, then you can modify them with hex flips to make them pretty (and get symmetry).
% By passing to a large enough edge inflation, it seems that we can always get n-fold symmetry, particularly by putting rhombs with angle 2 pi / n in the corners.  
% Symmetrize the KSK algorithm to make the result prettier.  
% This requires modifying the edges a bit so they only stick out, not in.
% Penrose boundary is special.
% Binary tiling is different.  

\section{Finding boundaries that work}\label{SEC:boundaries}

\subsection{Permutations}\label{SUBSEC:permutations}

It can be difficult to find an edge sequence $\seq$ that produces distorted boundaries $\seq\Rhomb{i}$ that all have tilable interiors.  
It would be good to have some general rules to help determine which sequences will work and which will not, but the approach used here is a simple brute-force computer search.  
This computer search begins with a standard edge sequence $\seq$ and permutes it, and for each permutation $\seq'$ applies the KSK criterion to see if all the $\seq'$-boundaries enclose tilable regions.  

\subsection{Iterators}\label{SUBSEC:iterators}

The edge sequence $\seq$ can be quite large---for instance, the elevenfold example in Section \ref{SEC:example} uses a sequence of length 35 with many repetitions---so finding all of its permutations involves non-trivial work.  
Specifically, it is not practical to store a list of all permutations of such a large sequence in memory; a much better idea is to use an iterator that can produce the next permutation---in some specified order---given the current permutation.  
Such an iterator only holds a single permutation in memory at one time, which reduces memory use enough to make it possible to search the permutations of very large sequences.  
The iterator used here is a variant of the cool-lex iterator \cite{RW:cool-lex} that has been adapted for iterating through permutations of multisets \cite{W:multisets}.  

\subsection{Simplifying the search}\label{SUBSEC:simplification}

Even when an iterator makes it feasible in theory to check the KSK condition for every permutation of an edge sequence, it may still take a prohibitively long time.  
For such cases, it is sometimes useful to reduce the number of permutations checked by restricting attention to permutations that are obtained by concatenating certain selected subsequences.  

For example, the elevenfold example in Section \ref{SEC:example} uses the edge sequence 
\begin{align*}
\seq & = (-1,1,-3,3,0,2,-2,-1,1,0,-5,5,-3,3,-1,1,4,-4,\\
& \qquad \qquad 2,-2,0,-1,1,2,-2,-3,3,0,4,-4,-1,1,2,-2,0).
\end{align*}

This edge sequence was found using a brute-force search as described above, but rather than checking all permutations of the following multiset 
\[
\begin{array}{rrrrrr}
0\times 5, &  1\times 5, &  2\times 4, &  3\times 3, &  4\times 2, &  5\times 1, \\
           & -1\times 5, & -2\times 4, & -3\times 3, & -4\times 2, & -5\times 1, 
\end{array}
\]
the search checks all sequences obtained by concatenating permutations of the following multiset of subsequences: 
\[
\begin{array}{rrr}
(0)\times 5, & (-1,1)\times 5, & (2,-2)\times 4, \\ (-3,3)\times 3, & (4,-4)\times 2, & (-5,5)\times 1.
\end{array}
\]
This drastically reduces the number of sequences that have to be checked.  

% An area calculation reveals that some boundaries will never work in any permutation.  
% Those sequences that are not rejected by this criterion need to be checked in all permutations.  
% You can apply the KSK algorithm by brute force.
% This will work in parallel: parcel out subsets of the list of all sequences to be searched by different threads.
% You can reduce your work by insisting that certain subsequences appear beside each other: -1 1, -2, 2, -3 3, etc..
% Substitution matrices are guaranteed to be symmetric.  

\section{An example with $11$-fold symmetry}\label{SEC:example}

This section is devoted to a single example of a substitution rule that yields tilings with elevenfold rotational symmetry.  
To my knowledge, it is the first such example published.  

The edge sequence is 
\begin{align*}
\seq & = (-1,1,-3,3,0,2,-2,-1,1,0,-5,5,-3,3,-1,1,4,-4,\\
& \qquad \qquad 2,-2,0,-1,1,2,-2,-3,3,0,4,-4,-1,1,2,-2,0).
\end{align*}
The associated inflation factor $\infl = \| \tot{\seq}\| \approx 27.2004$ is a Pisot number, which is a necessary condition for the existence of any non-trivial eigenvalue of the tiling dynamical system, that is, of any non-trivial discrete part in the spectrum (see \cite[Chapters 6 and 7]{BG:book} and \cite{S:pisot}).  

\begin{figure}
\includegraphics[width=\textwidth]{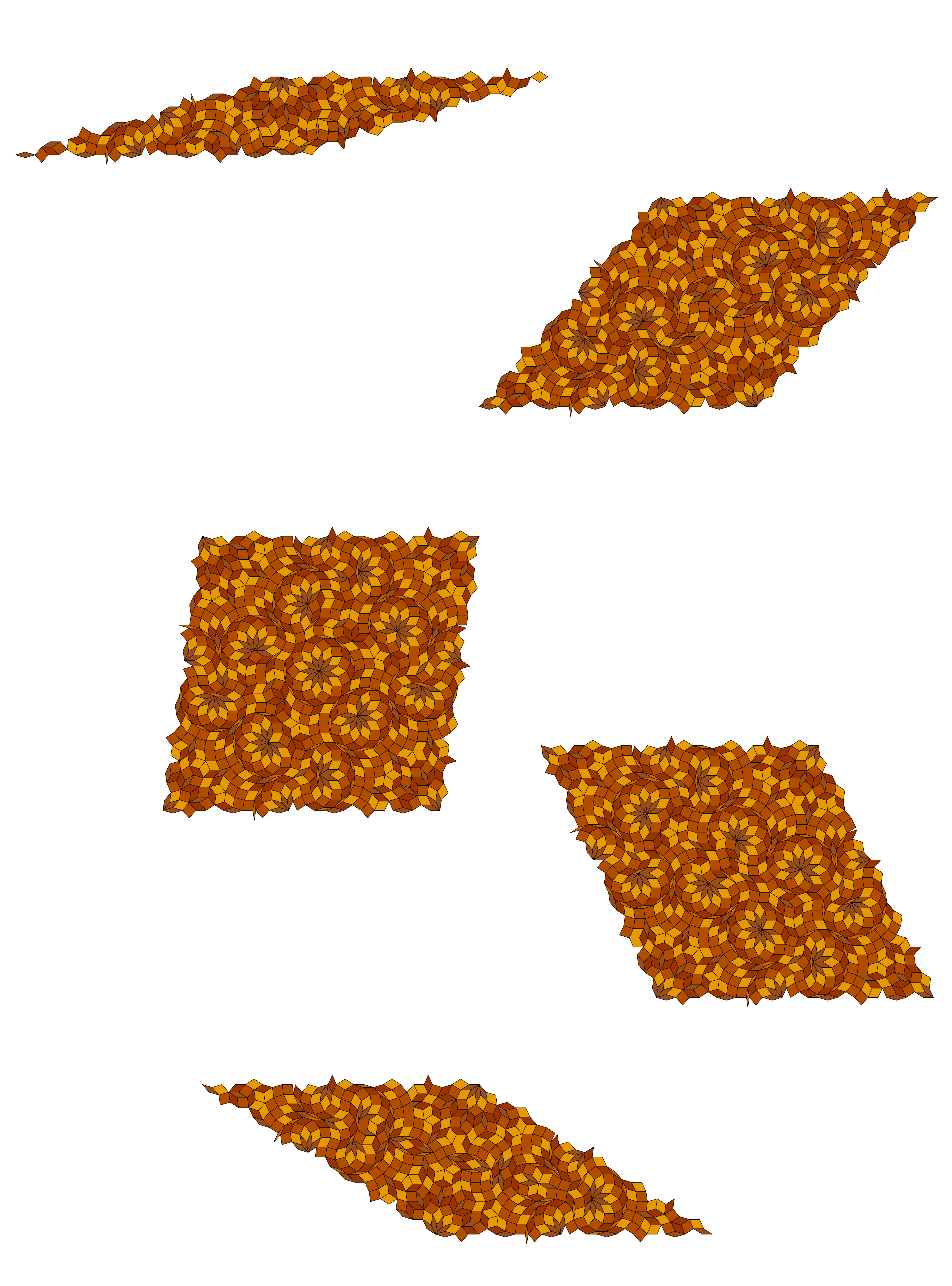}
\caption{\label{FIG:example-prototiles}Substituted images of $\Rhombs{11}$.  From top to bottom: $\sub (\Rhomb{10})$, $\sub(\Rhomb{8})$, $\sub(\Rhomb{6})$, $\sub(\Rhomb{4})$, and $\sub(\Rhomb{2})$.  }
\end{figure}

\begin{figure}
\includegraphics[width=\textwidth]{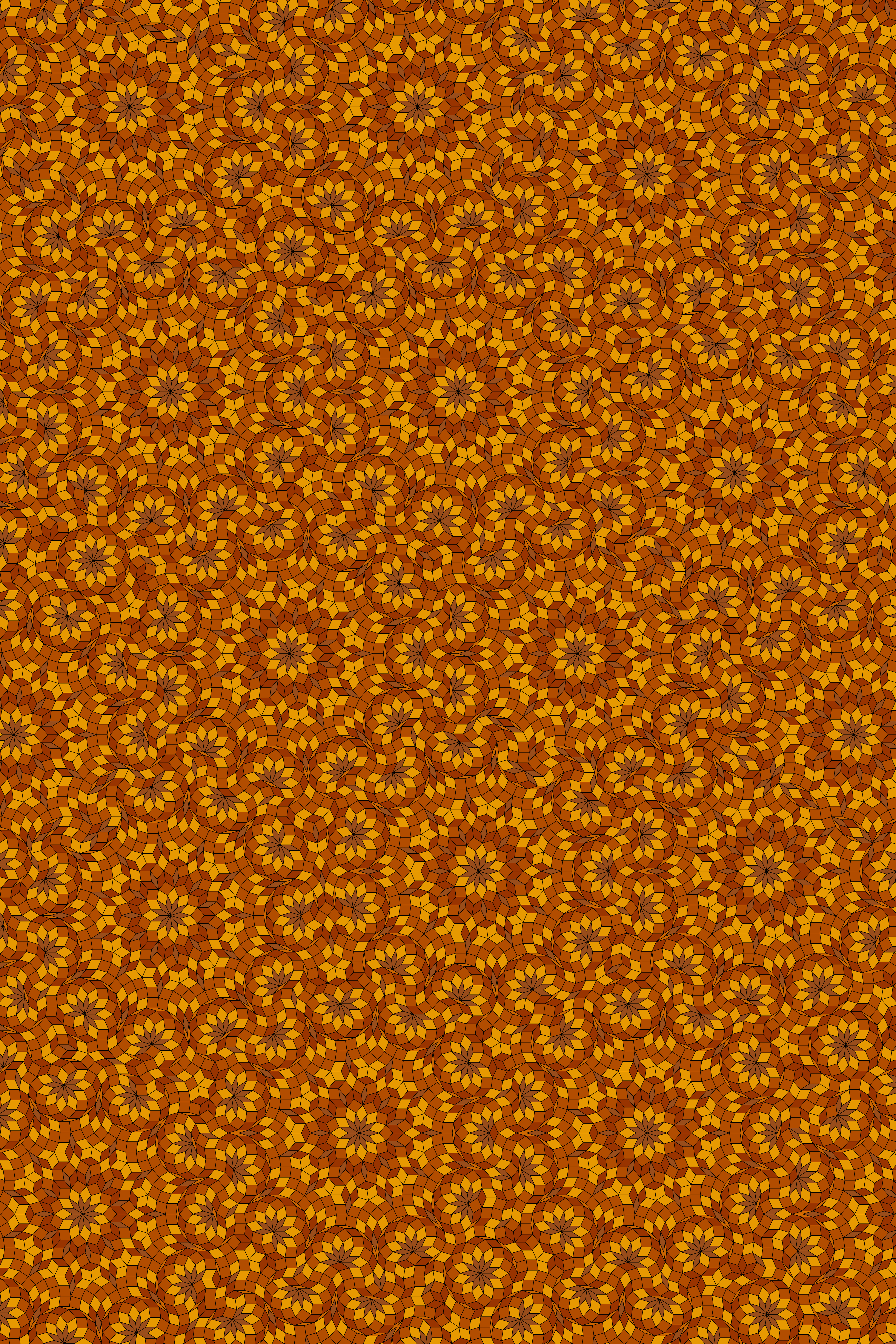}
\caption{\label{FIG:example-patch}A patch generated using the substitution rule in Figure \ref{FIG:example-prototiles}.}
\end{figure}

The substitution rule appears in Figure \ref{FIG:example-prototiles}, and a patch appears in Figure \ref{FIG:example-patch}.  

% Just some pictures. 

\section{The case when $n$ is even}\label{SEC:even}

The same methods can be applied to even $n$, but more care is necessary.  
Rule \ref{RULE} no longer applies as it did for odd $n$; indeed, if $n$ is even then a line segment parallel to $\edge_0$ is also parallel to $\edge_n$, and these two vectors would give it two opposite orientations.  
This means it is not so easy to determine the isometries that carry tiles in a given patch onto the prototiles.  

Therefore even after finding a tiling of the regions inside the $\seq$-boundaries of the rhombs $\Rhombs{n}$, one must still do more work to determine how the prototiles have been placed to obtain those tiles.  
I know of no way to automate this choice, but there are various ad hoc methods that have been successful \cite{B:ammann-beenker, H:rhomb}.  
One such method is to use prototiles in which opposite edges do not have the same orientation.  
This is what was done in \cite{H:rhomb} and in the version of the Ammann-Beenker tilings that appears in the Tilings Encyclopedia \cite{TE}.  

% How do the edge orientations on the prototiles work when n is even?

\bibliographystyle{abbrv}
\bibliography{bib-rhomb-tilings}

\begin{thebibliography}{10}

\bibitem{BG:book}
M.~Baake and U.~Grimm.
\newblock {\em Aperiodic order. {V}ol. 1}, volume 149 of {\em Encyclopedia of
  Mathematics and its Applications}.
\newblock Cambridge University Press, Cambridge, 2013.
\newblock A mathematical invitation, With a foreword by Roger Penrose.

\bibitem{BKSZ:projection}
M.~Baake, P.~Kramer, M.~Schlottmann, and D.~Zeidler.
\newblock Planar patterns with fivefold symmetry as sections of periodic
  structures in {$4$}-space.
\newblock {\em Internat. J. Modern Phys. B}, 4(15-16):2217--2268, 1990.

\bibitem{B:ammann-beenker}
F.~P.~M. Beenker.
\newblock Algebraic theory of non-periodic tilings of the plane by two simple
  building blocks: a square and a rhombus.
\newblock TH-Report 82-WSK04, Eindhoven University of Technology, the
  Netherlands, 1982.

\bibitem{F:thesis}
D.~Frettl\"{o}h.
\newblock Inflation\"{a}re {P}flasterungen der {E}bene mit
  ${D}_{2m+1}$-{S}ymmetrie und minimaler {M}usterfamilie.
\newblock Diploma {T}hesis, Universit\"{a}t Dortmund, 1998.

\bibitem{FH:worms}
D.~Frettl{\"o}h and E.~Harriss.
\newblock Parallelogram tilings, worms, and finite orientations.
\newblock {\em Discrete Comput. Geom.}, 49(3):531--539, 2013.

\bibitem{TE}
D.~Frettl\"{o}h and E.~O. Harriss.
\newblock The tilings encyclopedia.
\newblock Accessed: 2014-04-01.

\bibitem{G:shield}
F.~G\"{a}hler.
\newblock Crystallography of dodecagonal quasicrystals.
\newblock In C.~Janot and J.~Dubois, editors, {\em Quasicrystalline materials}.
  Institut Laue-Langevin, Committee for the Development in Europe of Science
  and Technology, World Scientific, 1988.
\newblock Proceedings of the I.L.L./CODEST workshop, Grenoble, 21-25 March
  1988.

\bibitem{GKM:search}
F.~G\"{a}hler, E.~Kwan, and G.~Maloney.
\newblock A computer search for planar substitution tilings with $n$-fold
  rotational symmetry.
\newblock {\em CoRR}, 2014.

\bibitem{G:multi}
J.~Garc{\'{\i}}a~Escudero.
\newblock Random tilings of spherical 3-manifolds.
\newblock {\em J. Geom. Phys.}, 58(11):1451--1464, 2008.

\bibitem{G:fivefold}
J.~Garc{\'{\i}}a~Escudero.
\newblock Randomness and topological invariants in pentagonal tiling spaces.
\newblock {\em Discrete Dyn. Nat. Soc.}, pages Art. ID 946913, 23, 2011.

\bibitem{H:rhomb}
E.~O. Harriss.
\newblock Non-periodic rhomb substitution tilings that admit order {$n$}
  rotational symmetry.
\newblock {\em Discrete Comput. Geom.}, 34(3):523--536, 2005.

\bibitem{KS:parallelograms}
S.~Kannan and D.~Soroker.
\newblock Tiling polygons with parallelograms.
\newblock {\em Discrete Comput. Geom.}, 7(2):175--188, Mar. 1992.

\bibitem{K:parallelograms}
R.~Kenyon.
\newblock Tiling a polygon with parallelograms.
\newblock {\em Algorithmica}, 9(4):382--397, 1993.

\bibitem{LB:binary}
F.~Lan{\c{c}}on and L.~Billard.
\newblock Two-dimensional system with a quasicrystalline ground state.
\newblock {\em J. Physique}, 49(2):249--256, 1988.

\bibitem{MSRRSB:lasers}
J.~Mikhael, M.~Schmiedeberg, S.~Rausch, J.~Roth, H.~Stark, and C.~Bechinger.
\newblock Proliferation of anomalous symmetries in colloidal monolayers
  subjected to quasiperiodic light fields.
\newblock {\em Proceedings of the National Academy of Sciences},
  107(16):7214--7218, 2010.

\bibitem{DN:n-fold}
K.-P. Nischke and L.~Danzer.
\newblock A construction of inflation rules based on {$n$}-fold symmetry.
\newblock {\em Discrete Comput. Geom.}, 15(2):221--236, 1996.

\bibitem{P:penrose}
R.~Penrose.
\newblock {The R\^ole of Aesthetics in Pure and Applied Mathematical Research}.
\newblock {\em The Institute of Mathematics and its Applications Bulletin},
  10(7/8):266--271, July 1974.

\bibitem{R:subrosa}
M.~Rissanen.
\newblock Sub rosa, a system of non-periodic rhombic substitution tilings.
\newblock Unpublished, 2013.

\bibitem{RW:cool-lex}
F.~Ruskey and A.~Williams.
\newblock The coolest way to generate combinations.
\newblock {\em Discrete Math.}, 309(17):5305--5320, 2009.

\bibitem{S:nobel}
D.~Shechtman, I.~Blech, D.~Gratias, and J.~W. Cahn.
\newblock Metallic phase with long-range orientational order and no
  translational symmetry.
\newblock {\em Phys. Rev. Lett.}, 53:1951--1953, Nov 1984.

\bibitem{S:socolar}
J.~E.~S. Socolar.
\newblock Simple octagonal and dodecagonal quasicrystals.
\newblock {\em Phys. Rev. B}, 39:10519--10551, May 1989.

\bibitem{S:pisot}
B.~Solomyak.
\newblock Dynamics of self-similar tilings.
\newblock {\em Ergodic Theory Dynam. Systems}, 17(3):695--738, 1997.

\bibitem{WSI:dodecagonal}
Y.~Watanabe, T.~Soma, and M.~Ito.
\newblock A new quasiperiodic tiling with dodecagonal symmetry.
\newblock {\em Acta Crystallographica Section A}, 51(6):936--942, Nov 1995.

\bibitem{W:multisets}
A.~Williams.
\newblock Loopless generation of multiset permutations using a constant number
  of variables by prefix shifts.
\newblock In {\em Proceedings of the Twentieth Annual ACM-SIAM Symposium on
  Discrete Algorithms}, SODA '09, pages 987--996, Philadelphia, PA, USA, 2009.
  Society for Industrial and Applied Mathematics.

\end{thebibliography}
\end{document}